\documentclass[]{interact}

\usepackage{epstopdf}% To incorporate .eps illustrations using PDFLaTeX, etc.
\usepackage[caption=false]{subfig}% Support for small, `sub' figures and tables

\usepackage[numbers,sort&compress]{natbib}% Citation support using natbib.sty
\bibpunct[, ]{[}{]}{,}{n}{,}{,}% Citation support using natbib.sty
\theoremstyle{plain}% Theorem-like structures provided by amsthm.sty
\newtheorem{theorem}{Theorem}[section]
\newtheorem{lemma}[theorem]{Lemma}

\newtheorem{proposition}[theorem]{Proposition}

\theoremstyle{definition}

\theoremstyle{remark}

\usepackage{mcode}
\usepackage{algorithm}
\usepackage[noend]{algpseudocode}
\algrenewcommand{\algorithmiccomment}[1]{{\hfill\color{black}$\blacktriangleright$ #1}}
\def\diag{\mathop{\mathrm{diag}}}
\def\max{\mathop{\mathrm{max}}}

\def\cond{\mathop{\mathrm{cond}}}

\def\vec{\mathop{\mathrm{vec}}}

\def\acos{\mathop{\mathrm{acos}}}
\begin{document}

%\articletype{ARTICLE TEMPLATE}% Specify the article type or omit as appropriate

\title{A Technique for Improving the Computation of Functions of Triangular Matrices}

\author{
\name{Jo\~ao R. Cardoso\textsuperscript{a}\thanks{Corresponding author: A. Sadeghi. Email: drsadeghi.iau@gmail.com} and Amir Sadeghi\textsuperscript{b}}
\affil{\textsuperscript{a}Coimbra Polytechnic/ISEC, Portugal, and Center for Mathematics, University of Coimbra, Portugal; \\
			\textsuperscript{b}Department of Mathematics, Parand and Robat Karim Branch, Islamic Azad University, Tehran, Iran.}
}

\maketitle

\begin{abstract}
	We propose a simple technique that, if combined with algorithms for computing functions of triangular matrices, can make them more efficient. Basically, such a technique consists in a specific scaling similarity transformation that reduces the departure from normality of a triangular matrix, thus decreasing its norm and in general its function condition number. It can easily be extended to non-triangular matrices, provided that it is combined with algorithms involving a prior Schur decomposition. Situations where the technique should be used or not will be discussed in detail. Special attention is devoted to particular algorithms like the inverse scaling and squaring to the matrix logarithm and the scaling and squaring to the matrix exponential. The advantages of our proposal are supported by theoretical results and illustrated with numerical experiments.
\end{abstract}

\begin{keywords}
Triangular matrices; Schur decomposition; matrix functions; scaling; condition number; matrix exponential; matrix logarithm; matrix square roots; matrix inverse cosine.
\end{keywords}

\section{Introduction}

Given a square complex matrix $A \in \mathbb{C}^{n \times n}$ and a scalar valued function $f$ defined on the spectrum of $A$ \cite[Def.1.1]{Higham}, the notation $f(A)$ means that {\it ``$f$ is a primary function of $A$''} in the usual sense, as considered in \cite{Higham} and \cite[Ch.6]{Horn}. We refer the reader to those books for background on matrix functions. Most of the methods for computing $f(A)$ require an initial Schur decomposition $A=UTU^*$ ($U$ unitary and $T$ upper triangular) so that the problem of computing $f(A)$ is then reduced to that of computing $f(T)$, because $f(A)=U\,f(T)\,U^*$.

We recall that many applications in science and engineering involve matrix functions computation; see \cite{Higham} for a large set of examples of applications and also the recent papers \cite{Arioli,Benzi,Rossignac} for additional applications.

In this paper, we begin by describing, in Section \ref{pre}, a scaling technique that can be easily incorporated into algorithms for computing functions of triangular matrices. As shown in Section \ref{sec-properties}, the main point is to reduce the norm of both $T$ and $f(T)$, and, in general, of the relative condition number of $f$ at $T$, thus making the problem of approximating $f(T)$ better conditioned. If conveniently implemented, this technique may bring benefits in terms of efficiency, specially in algorithms involving scaling and squaring (e.g., matrix exponential or direct trigonometric/hyperbolic matrix functions) or inverse scaling and squaring techniques (matrix logarithm, matrix $p$th roots and  inverse trigonometric/hyperbolic matrix functions).
Situations where the proposed technique has or does not have particular advantages are discussed and seem to be related with the size of the ratio $\|N\|_F/\|D\|_F$, where $D$ and $N$, are, respectively, the diagonal and the nilpotent parts of $T$.	
 In Section \ref{issues}, we discuss many issues related with the implementation of the scaling technique in finite precision environments, with particular emphasis on the computational cost (of order $O(n^2)$), error analysis and the proposal of a practical algorithm suitable to be combined with algorithms for matrix functions computation, which provides a heuristic to choose the scaling parameter $\alpha$. A set of experiments regarding the scaling and squaring method for the matrix exponential, the inverse scaling and squaring to the matrix logarithm and an algorithm for computing the inverse cosine of a matrix are carried out in Section \ref{experiments}. In Section \ref{conclusions}, some conclusions are drawn.
Unless otherwise stated, $\|.\|$, $\|.\|_F$ and $\|.\|_2$ will denote, respectively, a general subordinate matrix norm, the Frobenius norm and the 2-norm (also known as spectral norm). For a given $A \in \mathbb{C}^{n \times n}$, $|A|$ stands for the matrix of the absolute values of the entries of $A$.

\section{The Scaling Technique}\label{pre}

Let $T\in\mathbb{C}^{n\times n}$ be an upper triangular matrix and let $\alpha$ be a positive scalar. Without loss of generality, let us assume throughout the paper that $\alpha>1$. Consider the diagonal matrix
\begin{equation}\label{S}
S=\diag(1,\alpha,\ldots,\alpha^{n-1})
\end{equation}
and let us denote $\widetilde{T}:=S\,T\,S^{-1}$.
\medskip Given a complex valued function $f$ defined on the spectrum of $T$, the following steps describe, at a glance, the scaling technique we are proposing for computing $f(T)$:

\begin{enumerate}
	\item Choose a suitable scalar $\alpha>1$;
	\item Compute $f(\widetilde{T})$ using a certain algorithm;
	\item Recover $f(T)=S^{-1}\, f(\widetilde{T})\, S$.
\end{enumerate}

A discussion on the choice of $\alpha$ will be provided in  Section \ref{practical}. Step 3 is based on the fact that primary matrix functions preserve similarity transformations (see \cite[Thm. 1.13]{Higham}).

Note that the similarity transformation in Step 3 can magnify the errors arising in the computation of $f(\widetilde{T})$, thus resulting in a larger error in the approximation of $f(T)$. This issue will be discussed in Section \ref{error}, where we see that such errors are not larger that the ones resulting from the direct application of the algorithms (i.e., without using the scaling technique).

\medskip To gain insight into the effects of the left multiplication of $T$ by $S$ and of the right multiplication by $S^{-1}$, we write
\begin{equation}\label{N1N2}
T=D+N_1+\cdots+N_{n-1},
\end{equation}
where $D$ is a diagonal matrix formed by the diagonal of $T$ and zeros elsewhere, $N_1$ is formed by the first super-diagonal of $T$ and zeros elsewhere and so on, up to $N_{n-1}$.
Then
$$\widetilde{T}=STS^{-1}=D+ N_1/\alpha+\cdots+ N_{n-1}/\alpha^{n-1},$$
which means that the proposed technique just involves multiplications/divisions of certain entries of the matrix by the positive scalar $\alpha$.

For instance, if  $T=[t_{ij}]$ ($i,j=1,2,3$) is an $3\times 3$ upper triangular matrix ($t_{ij}=0$ for $i>j$), we have
$$D=\left[\begin{array}{ccc}
	t_{11} & 0 & 0\\
	0 & t_{22} & 0\\
	0 & 0 & t_{33}
\end{array}\right], \quad
N_1=\left[\begin{array}{ccc}
	 0 & t_{12} & 0\\
	0 & 0 & t_{23} \\
	0 & 0 &0
\end{array}\right], \quad
N_2=\left[\begin{array}{ccc}
0 & 0 & t_{13}\\
0 & 0 & 0\\
0 & 0 & 0
\end{array}\right].$$
Hence,
$$ \widetilde{T}=D+ N_1/\alpha+N_2/\alpha^2=
\left[\begin{array}{ccc}
t_{11} & t_{12}/\alpha & t_{13}/\alpha^2\\
0 & t_{22} & t_{23}/\alpha \\
0 & 0 &t_{33}
\end{array}\right].$$

\medskip To illustrate the benefits that the proposed scaling technique may bring, let us consider the problem of computing the exponential of the matrix
\begin{equation}\label{T}
T=\left[\begin{array}{cc}
1 & 10^6 \\
0 & -1 \\
\end{array}\right]
\end{equation}
by the classical scaling and squaring method, as described in \cite{Moler}. Before using Taylor or Pad\'e approximants, $T$ has to be scaled by $2^{k_0}$ so that the condition $\|T\|_2/2^{k_0}<1$ holds. We easily find  that the smallest $k_0$ verifying that condition is $k_0=21$, thus making it necessary to carry out at least 21 squarings for evaluating $e^T$. In contrast, if we apply the proposed  technique with $\alpha=10^{6}$, it is enough to take $k_0=1$, meaning that the computation of $e^T$ involves only one squaring and the multiplication of the $(1,2)$ entry of $\widetilde{F}=e^{\widetilde{T}}$ by $\alpha=10^6$, which is a very inexpensive procedure.

In the last decades, the scaling and squaring method has been subject to significant improvements; see, in particular, \cite{Mohy09,Higham05}. The function \texttt{expm} of the recent versions of MATLAB implements the method proposed in \cite{Mohy09}, where a sophisticated technique is used to find a suitable number of squarings. Such a technique is based on the magnitude of $\|T^k\|^{1/k}$ instead of $\|T\|$, which may lead to a considerable reduction in the number of squarings. For instance, if we compute $e^T$, where $T$ is the matrix given in (\ref{T}), by \cite[Alg. 5.1]{Mohy09}, no squaring is required. Note that this does not represent a failure of our proposal, because the scaling technique described above can be combined easily with any method for approximating the matrix exponential (or any other primary matrix function), in particular, with the new scaling and squaring method of \cite{Mohy09}. For instance, the computation of the exponential of the matrix in (\ref{T1}) involves $4$ squarings if computed directly by \texttt{expm} (which implements \cite[Alg. 5.1]{Mohy09}) and no squaring if scaled previously with $\alpha=3\times 10^4$. The main reason for saving squarings is that for $\alpha>1$, we have  $\|\widetilde{T}\|\leq \|T\|$ and consequently,
 $$\|\widetilde{T}^k\|^{1/k} \leq \|T^k\|^{1/k},$$
 for any positive integer $k$.
Similarly, our technique can be used to reduce the number of square roots  in the inverse scaling and squaring method to compute the matrix logarithm. The function \texttt{logm} implements the method provided in \cite[Alg. 5.2]{Mohy12}, where, like the matrix exponential, the estimation of the number of square roots is based on the magnitude of $\|T^k\|^{1/k}$ instead of $\|T\|$. To illustrate the gain of scaling a matrix, let us consider
 \begin{equation}\label{T1}
 T=\left[\begin{array}{cccc}
\mathtt{3.2346e-001} & \mathtt{3.0000e+004} & \mathtt{3.0000e+004} & \mathtt{3.0000e+004}\\
 0 & \mathtt{3.0089e-001} & \mathtt{3.0000e+004} & \mathtt{3.0000e+004}\\
 0 & 0 & \mathtt{3.2210e-001} & \mathtt{3.0000e+004}\\
 0 & 0 & 0 & \mathtt{3.0744e-001}
\end{array}\right]
\end{equation}
(see \cite[p.C163]{Mohy12}). The direct computation of $\log(T)$ by \texttt{logm} requires the computation of $16$ square roots, while an initial scaling of $T$ with $\alpha=3\times 10^4$, requires only $5$ square roots, without any sacrifice in the accuracy of the computed logarithm (see the results for $T_2$ in Figure \ref{fig-log}).

We finish this section by noticing that this scaling technique may also be combined with many other methods that reduce a matrix to upper triangular form in order to perform the computation: the Schur-Parlett method of \cite{Davies} for computing several matrix functions (available through the MATLAB function \texttt{funm}) and the algorithms  provided in \cite{Bjorck} for the matrix square root, in \cite{Higham13,Iannazzo,Smith} for $p$th roots or fractional powers of a matrix, in \cite{Aprahamian14} for the unwinding function, in \cite{Mohy15} for the matrix sine and cosine, in \cite{Fasi} for the Lambert W function and in \cite{Cardoso,Garrappa} for special matrix functions.

\section{Properties}\label{sec-properties}

To understand why the proposed scaling technique may improve the matrix function computations, we show, in this section, that, for any $\alpha>1$, it reduces the Frobenius norms of both $T$ and of $f(T)$. We also provide insight to understand why in all the examples we have tested the norm of the Fr\'echet derivative $L_f(T)$ is reduced as well, thus making the problems better conditioned.

\begin{proposition}\label{property1}
	 Let us assume that: $T\in\mathbb{C}^{n\times n}$ is an upper triangular matrix, $\alpha>1$, $f$ is a complex valued function defined on the spectrum of $T$, and  $S$ is defined as in (\ref{S}). Denoting $\widetilde{T}=STS^{-1}$, the following inequalities hold:
	 \begin{enumerate}
	 	\item[(i)] $\|\widetilde{T}\|_F\leq \|T\|_F$;
	 	\item[(ii)] $\|f(\widetilde{T})\|_F\leq \|f(T)\|_F$.	 	
	 \end{enumerate}
 \end{proposition}
	
\begin{proof}
\begin{enumerate}
	\item[(i)]
	Let us write $T=D+N_1+\cdots+N_{n-1}$, where $N_j,\ j=1,\ldots,n-1$, are defined as in (\ref{N1N2}). Then
	\begin{equation}\label{nT2}
	\|T\|_F^2=\|D\|_F^2+ \|N_1\|_F^2+\cdots+\|N_{n-1}\|_F^2.
	\end{equation}
	Since
	\begin{eqnarray*}
	\widetilde{T} &=& STS^{-1}\\
				&=& D+SN_1S^{-1}+\cdots+SN_{n-1}S^{-1}\\
				&=& D+N_1/\alpha+\cdots+N_{n-1}/\alpha^{n-1},
	\end{eqnarray*}
		we have
			\begin{equation}\label{n-tildeT2}
		\|\widetilde{T}\|_F^2=\|D\|_F^2+ \frac{1}{\alpha^2}\|N_1\|_F^2+\cdots+\frac{1}{\alpha^{2(n-1)}}\|N_{n-1}\|_F^2.
		\end{equation}	
	From (\ref{nT2}) and (\ref{n-tildeT2}),
	\begin{equation}\label{T-Ttilde}
	\|T\|_F^2=\|\widetilde{T}\|_F^2+ \left(1-\frac{1}{\alpha^2}\right)\|N_1\|_F^2+\cdots+\left(1-\frac{1}{\alpha^{2(n-1)}}\right)\|N_{n-1}\|_F^2.
	\end{equation}
	Since $\alpha>1$, the inequality $\|\widetilde{T}\|_F \leq \|T\|_F$ follows.

\item[(ii)]  Let us write $f(T)=f(D)+F_1+\cdots+F_{n-1},$
where $F_1$ is formed by the first super-diagonal of $f(T)$ and zeros elsewhere and so on, up to $F_{n-1}$.
Then
\begin{equation}\label{nfT2}
\|f(T)\|_F^2=\|f(D)\|_F^2+ \|F_1\|_F^2+\cdots+\|F_{n-1}\|_F^2,
\end{equation}
and
\begin{equation}\label{n-ftildeT2}
\|f(\widetilde{T})\|_F^2=\|f(D)\|_F^2+ \frac{1}{\alpha^2}\|F_1\|_F^2+\cdots+\frac{1}{\alpha^{2(n-1)}}\|F_{n-1}\|_F^2.
\end{equation}	
From (\ref{nfT2}) and (\ref{n-ftildeT2}),
\begin{equation}\label{T-fTtilde}
\|f(T)\|_F^2=\|f(\widetilde{T})\|_F^2+ \left(1-\frac{1}{\alpha^2}\right)\|F_1\|_F^2+\cdots+\left(1-\frac{1}{\alpha^{2(n-1)}}\right)\|F_{n-1}\|_F^2.
\end{equation}
Since $\alpha>1$, the result follows.	
\end{enumerate}   
\end{proof}
	
Using (\ref{T-Ttilde}), we can evaluate the difference between the squares of the norms $\|\widetilde{T}\|_F$ and $\|T\|_F$. Likewise, the difference between the squares of $\|f(\widetilde{T})\|_F$ and $\|f(T)\|_F$ can be found from (\ref{T-fTtilde}). In particular, if $T$ is nonsingular, we have $\|\widetilde{T}^{-1}\|_F\leq \|T^{-1}\|_F$,
which shows that $\kappa(\widetilde{T}) \leq \kappa(T)$, where $\kappa(X)=\|X\|\|X^{-1}\|$ denotes the standard condition number of $X$ with respect to the matrix inversion. This issue has also motivated us to investigate the absolute and relative condition numbers of $f$ at $\widetilde{T}$, which will be defined below.

Given a map $f:\mathbb{C}^{n\times n}\rightarrow\mathbb{C}^{n\times n}$, the Fr\'{e}chet derivative of $f$ at $A\in\mathbb{C}^{n\times n}$ in the direction of $E\in\mathbb{C}^{n\times n}$ is a linear operator $L_f(A)$ that maps the ``direction matrix'' $E$ to $L_f(A,E)$ such that
$$\lim_{E\rightarrow 0}\frac{\|f(A+E)-f(A)-L_f(A,E)\|}{\|E\|}=0.$$ The Fr\'{e}chet derivative of the matrix function $f$ may not exist at $A$, but if it does it is unique and coincides with the directional (or G\^ateaux) derivative of $f$ at $A$ in the direction $E$. Hence, the existence of the Fr\'echet derivative guarantees that, for any $E\in\mathbb{C}^{n\times n}$,
$$L_f(A,E)=\lim_{t\rightarrow 0}\frac{f(A+tE)-f(A)}{t}.$$
Any consistent matrix norm $\|.\|$ on $\mathbb{C}^{n\times n}$ induces the operator norm $\ \|L_f(A)\|:=$ \linebreak $\max_{\|E\|=1}\,\|L_f(A,E)\|.$
Here we use the same notation to denote both the matrix norm and the induced operator norm. The Fr\'echet derivative plays a key role to understand how $f(A)$ changes as $A$ is being subject to perturbations of first order. Its norm is commonly used to define the absolute and the relative condition numbers of $f$ at $A$:
\begin{eqnarray}
	\cond\,\hspace*{-0.5ex}_\mathrm{abs}(f,A)&=&\|L_f(A)\|;\label{cond-abs}\\
	\cond\,\hspace*{-0.5ex}_\mathrm{rel}(f,A)&=&\|L_f(A)\|\frac{\|A\|}{\|f(A)\|}.\label{cond-rel}
\end{eqnarray}

Since $L_f(A,E)$ is linear in $E$, it is often important to consider the so-called Kronecker form of the Fr\'echet derivative:
\begin{equation}\label{frechet-vec}
\vec\left(L_{f}(A,E)\right)=K_{f}(A)\,\vec(E),
\end{equation}
where $\vec(.)$ stands for the operator that stacks the columns of a matrix into a long vector of size $n^2\times 1$, and $K_{f}(A)\in \mathbb{C}^{n^2\times n^2}.$ With respect to the Frobenius norm, the following equality is valid:
\begin{equation}\label{equality-frechet}
	\|L_f(A)\|_F = \|K_f(A)\|_2.
\end{equation}

For details on the Fr\'echet derivative and its properties see, for instance, \cite[Ch. 3]{Higham} and also \cite[Ch. X]{Bhatia97}.

\begin{lemma}\label{lemma1}
	For each $i,j=1,\ldots,n$, let $E_{ij}$ be the matrix of order $n$ with $1$ in the position $(i,j)$ and zeros elsewhere. Using the same notation of Proposition \ref{property1}, the following equality holds:
	\begin{equation}\label{equal-frechet}
		L_f(\widetilde{T},E_{ij})=\alpha^{j-i}\,S\, L_f(T,E_{ij})\,S^{-1}.
	\end{equation}
\end{lemma}

\begin{proof}
	Through a simple calculation, we can see that, for any $i,j=1,\ldots,n$, $S^{-1}E_{ij}S=\alpha^{j-i} E_{ij}.$
	Now, by the linearity and similarity properties of Fr\'echet derivatives, one arrives at:
	\begin{eqnarray*}
		L_f(\widetilde{T},E_{ij})&=& L_f(STS^{-1},E_{ij})\\
		&=& S\, L_f(T,S^{-1}E_{ij}S)\,S^{-1}\\
		&=& \alpha^{j-i}\,S\, L_f(T,E_{ij})\,S^{-1}.
	\end{eqnarray*}
	
\end{proof}

Given any matrix $X$ of order $n$ (not necessarily triangular), let us consider the  following decomposition:
\begin{equation}\label{dec-X1}
	X = N_{-(n-1)}+\cdots + N_{-1}+N_0+N_1+\cdots+N_{n-1},
\end{equation}
where $N_0$ denote the matrix of order $n$ formed by the diagonal of $X$ and zeros elsewhere, and, for $p=1,\ldots,n-1$, $N_{-p}$ (resp., $N_p$) is formed by the $p$-th sub-diagonal of $X$ and zeros elsewhere (resp., the $p$-th super-diagonal of $X$ and zeros elsewhere).

 \begin{lemma}\label{lemma0}
 	For $i,j=1,\ldots,n$, let $E_{ij}$ be the matrix used in Lemma \ref{lemma1}.
 	If $T$ and $\widehat{T}$ are both upper triangular matrices of order $n$, let us denote $X:=T\,E_{ij}\,\widehat{T}$. Then, in the decomposition \eqref{dec-X1} of $X$, we have
 	$$N_{-(n-1)}=\cdots=N_{-(i-j+1)}=0,$$
 	that is,
 	\begin{equation}\label{dec-X}
 		X = N_{-(i-j)}+\cdots + N_{-1}+N_0+N_1+\cdots+N_{n-1}.
 	\end{equation}
 \end{lemma}

 \begin{proof}
 	We know that $E_{ij}={\mathbf e}_i{\mathbf e}_j^T$, where ${\mathbf e}_k$ is the $n\times 1$ vector with $1$ in the $k$-th position and zeros elsewhere.
 	Let ${\mathbf 0}_{p\times q}$ be the zero matrix of size $p\times q$. Denoting by $T(:,i)$ and $\widehat{T}(j,:)$, respectively, the $i$-th column of $T$ and the $j$-th row of $\widehat{T}$, some calculations show that
 	\begin{eqnarray*}
 		X &=& T\,E_{ij}\,\widehat{T} \\
 		& = & T{\mathbf e}_i{\mathbf e}_j^T \,\widehat{T} \\
 		& = & T(:,i)\,\widehat{T}(j,:) \\
 		& = & \left[\begin{array}{c}
 			t_{1i}\\
 			\vdots \\
 			t_{ii}\\
 			0 \\
 			\vdots \\
 			0
 		\end{array} \right]
 		\left[\begin{array}{cccccc}
 			0 & \cdots & 0 & \widehat{t}_{jj} & \cdots & \widehat{t}_{jn}
 		\end{array} \right] \\
 	& = & 	\left[	
 			\begin{array}{l|l}
 				{\mathbf 0}_{i\times (j-1)}  & P_{i\times (n-j+1)}\\ \hline
 				{\mathbf 0}_{(n-i)\times (j-1)} &{\mathbf 0}_{(n-i)\times (n-j+1)}		
 			\end{array}
 			\right],
 	 	\end{eqnarray*}
  	where $P$ is a matrix of size $i\times (n-j+1)$, from which the result follows.
  	
 \end{proof}

\begin{proposition}\label{property2}
	 Let us assume that: $T\in\mathbb{C}^{n\times n}$ is an upper triangular matrix, $\alpha>1$, $f$ is a complex valued function defined on the spectrum of $T$, and  $S$ is defined as in (\ref{S}). Denote $\widetilde{T}:=STS^{-1}$ and consider the matrix $E_{ij}$ defined in Lemma \ref{lemma1}. If, for any $i,j=1,\ldots,n$, there exists a matrix $P$ of size $i\times (n-j+1)$ such that
	 \begin{equation}\label{block-P}
	 L_f(T,E_{ij}) = \left[	
	 \begin{array}{l|l}
	 	{\mathbf 0}_{i\times (j-1)}  & P_{i\times (n-j+1)}\\ \hline
	 	{\mathbf 0}_{(n-i)\times (j-1)} &{\mathbf 0}_{(n-i)\times (n-j+1)}		
	 \end{array}
	 \right]
	 \end{equation}
	 then
	 	  \begin{equation}\label{inequal-frechet-2}
	 	 	\|L_f(\widetilde{T},E_{ij})\|_F\leq \|L_f(T,E_{ij})\|_F.
	 	 \end{equation}
	 	
\end{proposition}

\begin{proof}
From \eqref{block-P} and using the notation of Lemma \ref{lemma0}, we can write 	
$$L_f(T,E_{ij})= N_{-(i-j)}+\cdots + N_{-1}+N_0+N_1+\cdots+N_{n-1}.$$
Hence
	\begin{eqnarray*}
	L_f(\widetilde{T},E_{ij})&=&\alpha^{j-i}\,S\, L_f(T,E_{ij})\,S^{-1}\\
	     & = & \alpha^{j-i} \left(\alpha^{i-j}N_{-(i-j)}+\cdots + \alpha
	  N_{-1}+N_0+\frac{N_1}{\alpha}+\cdots+\frac{N_{n-1}}{\alpha^{n-1}}\right) \\
& = & N_{-(i-j)}+\frac{N_{-(i-j-1)}}{\alpha}\cdots + \frac{
N_{j-i+1}}{\alpha^{i-j-1}}+\frac{N_0}{\alpha^{i-j}}+\frac{N_1}{\alpha^{i-j+1}}+\cdots+\frac{N_{n-1}}{\alpha^{i-j+n-1}},
\end{eqnarray*}
from which the result follows.
  
\end{proof}
	
The Fr\'echet derivatives of the most used primary matrix functions are sums or integrals involving functions of the form $g(A,E)=f_1(A)\,E\,f_2(A)$, where $f_k(A)$ ($k=1,2$) is a certain primary matrix function. For instance, the Fr\'echet derivatives of the matrix exponential and matrix logarithm allow, respectively, the integral representations
\begin{equation}\label{frechet-exp}
L_{\exp}(A,E)=\int_{0}^1 e^{A(1-t)}Ee^{At}\ dt
\end{equation}
and
\begin{equation}\label{frechet-log}
L_{\log}(A,E)=\int_0^1\, \left(t(A-I)+I\right)^{-1}\,E\, \left(t(A-I)+I\right)^{-1}\ dt
\end{equation}
(see \cite{Higham}). More generally, a function that can be represented by a Taylor series expansion
$$
f(A)=\sum_{k=0}^\infty a_kA^k,
$$
has a Fr\'echet derivative of the form \cite{Kenney}
\begin{equation}\label{frechet-general}
L_f(A,E)=\sum_{k=0}^\infty a_k \sum_{j=0}^{k-1} A^{j}E A^{k-1-j},
\end{equation}
which involves sums of functions like $g(A,E)$.

 If $A$ is upper triangular, then, by Lemma \ref{lemma0} and Proposition \ref{property2}, $X=g(A,E_{ij})$ has a structure like \eqref{block-P}, which is preserved by sums or integrals. Hence, the Fr\'echet derivatives of the most known primary matrix functions including, in particular, (\ref{frechet-exp}), (\ref{frechet-log}) and (\ref{frechet-general}), verify the conditions of Proposition \ref{property2} and consequently the inequality \eqref{inequal-frechet-2}.

Let us now recall how to evaluate $\|L_f(A)\|_F$, which appears in the definition of the condition numbers of a function; check \eqref{cond-abs} and \eqref{cond-rel}. Once we know $L_f(A,E_{ij})$, for a given a pair $(i,j)$, with $i,j\in\{1,\ldots,n\}$, the equality (\ref{frechet-vec}) enables us to find the $((j-1)n+i)$-th column of $K_f(A)$. Repeating the process for all $i,j=1,\ldots,n$, we can find all the entries of $K_f(A)$. By \eqref{equality-frechet}, the absolute and relative conditions numbers follow easily.

We now compare the condition numbers corresponding to $\widetilde{T}$ and $T$, by analysing the values of $\|K_f(\widetilde{T})\|_2$ and $\|K_f(T)\|_2$. To simplify the exposition, let us denote $\widetilde{K}:=K_f(\widetilde{T})=[\widetilde{k}_{pq}]$ and  $K:=K_f(T)=[k_{pq}]$ ($p,q=1,\ldots,n^2$). Attending to (\ref{inequal-frechet-2}), the $2$-norm of the $p$-th column of $\widetilde{K}$ is smaller than or equal to that of the $p$-th column of $K$, for any $p=1,2,\ldots,n^2$. This means that,
$$\|\widetilde{K}{\mathbf e}_p\|_2 \leq \|K{\mathbf e}_p\|_2,$$
where ${\mathbf e}_p$ denotes the $n^2\times 1$ vector with one in the $p$-th component and zeros elsewhere. Moreover, $\diag(\widetilde{K})=\diag(K)$ and, by applying the $\vec$ operator to both hand-sides of (\ref{equal-frechet}), we can observe that $\widetilde{k}_{pq}$ and $k_{pq}$ have the same signs for all $p,q=1,2,\ldots,n^2$. If $\widetilde{K}$ (or $K$) has all the entries non-negative, the properties of the spectral norm ensure that
\begin{equation}\label{ineq-positive} 
	\|\widetilde{K}\|_2 \leq \|K\|_2,
\end{equation}
and, consequently, 
\begin{equation}\label{ineq-abs} 
	\cond\,\hspace*{-0.5ex}_\mathrm{abs}(f,\widetilde{T}) \leq \cond\,\hspace*{-0.5ex}_\mathrm{abs}(f,T).
\end{equation}
However, because the spectral norm is not absolute (that is, in general, $\|A\|_2\neq \|\ |A|\ \|_2$; see, for instance, \cite{Mathias} and \cite[Sec. 5.6]{Horn13}), in the case of $\widetilde{K}$ (or $K$) having simultaneously positive and negative entries, it is not fully guaranteed that (\ref{ineq-abs}) holds. Nevertheless, we believe that (\ref{ineq-abs}) is not valid just for a few exceptional cases. In all the tests we have carried out, we have not encountered any example for which (\ref{ineq-abs}) does not hold. The same can be said about the inequality 
$
\cond\,\hspace*{-0.5ex}_\mathrm{rel}(f,\widetilde{T}) \leq \cond\,\hspace*{-0.5ex}_\mathrm{rel}(f,T),
$
that our tests suggested to be true in general. 

\section{Implementation Issues}\label{issues}

\subsection{Computational Cost}

Provided that an algorithm (let us call it {\it Algorithm $x$}) for computing $f(T)$, where $f$ is a certain primary matrix function, $T$ is a given upper triangular matrix and $\alpha>1$ is a suitable scalar, the proposed technique is very easy to implement. Indeed, using the MATLAB code displayed on Figure \ref{fig1}, it can be carried out through the following steps:
\begin{enumerate}
	\item \texttt{T1 = multiply\_by\_alpha(T,alpha)}, for computing $\mathtt{T1}=STS^{-1}$.
	\item Run Algorithm $x$ to approximate $f(\mathtt{T1})$;
	\item Recover $F=f(T)\approx S^{-1}f(\mathtt{T1})S$ using $\ $\texttt{F = multiply\_by\_alpha(T1,1/alpha)}.
\end{enumerate}

\begin{figure}[ht]
	\centering
\begin{lstlisting}
function T1 = multiply_by_alpha(T,alpha)
n = length(T);
T1 = T;
for k = 1:n
    for s = k+1:n
        T1(k,s) = T1(k,s)*alpha^(s-k);
    end
end
end
\end{lstlisting}
	\caption{\small MATLAB code for performing a similarity transformation of an upper triangular matrix via $S$ in \eqref{S}.}
	\label{fig1}
\end{figure}

Note that the same code in Figure \ref{fig1} can be used to perform the two similarity transformations required by the technique.

Most of the effective algorithms for matrix functions involve $O(n^3)$ operations, while the proposed technique involves only $O(n^2)$. Thus, if the choice of $\alpha$ is such that some $O(n^3)$ operations are saved (e.g., squarings, square roots, matrix products,...), such a preconditioning technique can contribute towards a reduction in the computational cost.
\subsection{Error Analysis}\label{error}

Let us denote  $\widetilde{F}=f(\widetilde{T})=[\widetilde{f}_{ij}]$ and $F=f(T)=[f_{ij}]$, where $i,j=1,\ldots,n$. We assume that $\widetilde{X}\approx \widetilde{F}$ is the approximation arising in the computation of $f(\widetilde{T})$ using a certain algorithm (Algorithm $x$) and that $X=S^{-1}\widetilde{X}S$ is the approximation to $f(T)$ that results from the multiplication of $\widetilde{X}$ by $S^{-1}$ on the left-hand side and by $S$ on the right-hand side. Recall that $S$ is defined by (\ref{S}), where $\alpha>1$. The entries of $\widetilde{X}$ (resp., $X$) are denoted by $\widetilde{x}_{ij}$ (resp., $x_{ij}$).

Provided that an $\alpha$ is suitable chosen, it is expected that $f$ is well-conditioned at $\widetilde{T}$ so that the approximation $\widetilde{X}\approx \widetilde{F}$ is highly accurate. However, we must take into account that the similarity transformation $S^{-1} \widetilde{X}S$ will magnify the errors of the approximation $\widetilde{X}\approx \widetilde{F}$, so the success of our technique relies mainly on finding a suitable $\alpha$ that leads to a magnification of the errors not large than the corresponding errors arising in the direct application of an Algorithm $x$ for computing $f(T)$. To get more insight, let us analyze the component-wise absolute errors.
Let $\widetilde{\mathcal{E}}:=\widetilde{F}-\widetilde{X}=[\widetilde{\varepsilon}_{ij}]$ and $\mathcal{E}:=F-X=[\varepsilon_{ij}].$ We have
\begin{eqnarray}
	\mathcal{E} &=& F-X\nonumber\\
	&=& S^{-1}\widetilde{F}S- S^{-1}\widetilde{X}S\nonumber\\
	&=& S^{-1}(\widetilde{F}-\widetilde{X})S\nonumber \\
	&=& S^{-1}\widetilde{\mathcal{E}}S \nonumber\\
	&=& \left[
	\begin{array}{ccccc}
		\widetilde{\varepsilon}_{11}&\alpha\,\widetilde{\varepsilon}_{12}&\alpha^2\, \widetilde{\varepsilon}_{13}&\cdots&\alpha^{n-1}\,\widetilde{\varepsilon}_{1n}\\
		0&\widetilde{\varepsilon}_{22}&\alpha\,\widetilde{\varepsilon}_{23}&\ddots&\vdots\\
		\vdots&\ddots&\ddots&\ddots&\alpha^2\,\widetilde{\varepsilon}_{n-2,n}\\
		&&&&\alpha \,\widetilde{\varepsilon}_{n-1,n}\\
		0&\cdots&&0&\widetilde{\varepsilon}_{nn}\label{wise-error}
	\end{array}
	\right],
\end{eqnarray}
where we can see that some entries of $\widetilde{\mathcal{E}}$ are magnified by the powers of $\alpha$, so that $\alpha$ and $n$ cannot be both very large to avoid big errors in the top-right entries of $X$. For instance, if $\widetilde{\varepsilon}_{1n}=10^{-20}$ and $\alpha^{n-1}=10^{20}$, we get and error of order $1$! This illustrates that the success of our strategy for improving the computation of matrix functions requires a careful choice of $\alpha$, that should be done in accordance with the order $n$.

In the next subsection we discuss the choice of $\alpha$ and also propose a strategy to deal with matrices of large order $n$.

\subsection{Practical Implementation of the Technique} \label{practical}
Up to now, we know that $\alpha$ must be enough large for having a reduction in the norm and to get a smaller condition number but $\alpha^{n-1}$ cannot be much large to avoid introducing large errors or underflow/overflow. So if $n$ is large, $\alpha$ has to be small and so our technique may be useless. To round this issue, we will carry out the similarity transformations with a slightly different matrix $S$. We will work instead with
\begin{equation}\label{S-block}
	S=\diag(I_{n_1},\alpha\, I_{n_2},\ldots,\alpha^{m-1}\, I_{n_{m}}),
\end{equation}
where $I_{n_k}$ is the identity matrix of order $n_k$ and $n_1+n_2+\cdots + n_m=n$, which is formed by $m$ blocks that are multiples of the identity matrix.
This also allows a reduction on the norm of $T$, but now we can control $m$ independently of the size of $T$, making our technique applicable to larger matrices.
Finding an $\alpha$ that would be optimal for all functions $f$ and for any matrix $T$ may be very difficult, because an $\alpha$ appropriate for $f(T)$ may not be for other function, so our strategy is to provide a heuristic that combines the values of $\alpha$ and $m$, in order make our strategy effective for computations using IEEE standard double precision arithmetic.

To simplify, we assume that in \eqref{S-block} the first $m-1$ blocks have equal
 size $n_1=\lfloor\frac{n}{m}\rfloor$ (here, $\lfloor.\rfloor$ denotes the floor of a number), that is, $n_1=\ldots=n_{m-1}=\lfloor\frac{n}{m}\rfloor$ and $n_m=n-(m-1)n_1$. To illustrate, for $n=57$ and $m=6$, we have
$$S=\diag(I_9,\alpha\, I_9,\alpha^2\, I_9,\alpha^3\, I_9,\alpha^4\, I_9,\alpha^5\, I_{12}).$$
To find appropriate values for $\alpha$ and $m$, our experience with a large set of tests carried out using IEEE standard double precision arithmetic suggests taking $\alpha$ as the maximal absolute value of $T$:
$$\alpha=\max \{|t_{ij}|:\ i=1,\ldots,n,\ j=i,\ldots,n \}$$
and $m$ as the largest integer such that
$$\alpha^m \leq 10^{20}.$$ To avoid situations where the scaling by $\alpha$ does not bring any benefits, we assume that $\alpha\geq 10$. Of course, other options for $\alpha$ and $m$ may be considered, depending on the function and on the matrix we are dealing with, and also on the precision used.
Algorithm \ref{alg1} summarizes a practical implementation of our technique.

\begin{algorithm}[t]
	\caption{This algorithm proposes a scaling technique for improving the performance of a certain Algorithm $x$ for computing $Y=f(T)$, where $T$ is an upper triangular matrix of order $n$. }
	\label{alg1}
	{
		\begin{algorithmic}[1]
			\State $\alpha \gets \max \{|t_{ij}|:\ i=1,\ldots,n,\ j=i,\ldots,n \}$;
			\State $m\gets \lfloor{20\log(10)/\log(\alpha)}\rfloor$;
			\State If $m>n$, set $m\gets n$;
			\State If $\alpha<10$, set $\alpha\gets 1$;
			\State $n_1\gets \lfloor{n/m}\rfloor$ and $n_m\gets n-(m-1)n_1$;
			\State $S = \diag(I_{n_1},\alpha\, I_{n_1},\ldots,\alpha^{m-2}\, I_{n_1}, \alpha^{m-1}\, I_{n_m});$
			\State $S^{-1} = \diag(I_{n_1},\alpha^{-1}\, I_{n_1},\ldots,\alpha^{2-m}\, I_{n_1}, \alpha^{1-m}\, I_{n_m});$
			\State $\widetilde{T}\gets STS^{-1}$;
			\State Compute $\widetilde{Y} \gets f(\widetilde{T})$ using Algorithm $x$;
			\State $Y\gets S^{-1}\widetilde{Y}S.$
			\end{algorithmic}
	}
\end{algorithm}

\section{Numerical experiments}\label{experiments}

We have implemented Algorithm \ref{alg1} in MATLAB R2021a (with unit roundoff $u\approx 1.1\times 10^{-16}$) and report below on the relative error and on the number of squarings/inverse squarings, for several experiments involving the computation of the exponential, the logarithm and the inverse cosine of upper triangular matrices. 

The following notation will be used:

\begin{itemize}
	\item[$\bullet$] $s$ and $\widetilde{s}$ denote, respectively, the number of squarings without and with the scaling technique required in the computation of the matrix exponential or the number of square roots involved in the inverse scaling and squaring method for the matrix logarithm and for the inverse cosine;
	\item[$\bullet$] $e$ and $\widetilde{e}$ are the relative errors of the computed approximations with respect to the Frobenius norm, that is,
	\begin{eqnarray*}
		e&=&\|f(T)-\widehat{f}(T)\|_F/\|f(T)\|_F\\
		\widetilde{e}&=&\|f(T)-S\,\widehat{f}(\widetilde{T})\,S^{-1}\|_F/\|f(T)\|_F,
	\end{eqnarray*}
where $f(T)$ means the exact function of $T$ and $\widehat{f}(T)$ is the approximation obtained by a certain Algorithm $x$. In all the experiments we also compute the values of the ratio $\|N\|_F/\|D\|_F$ for illustrating that savings in the numbers of square roots or squarings are expected when this ratio is large, suggesting that the success of our technique might be related with those values.

\end{itemize}

\noindent {\it Experiment 1.} In this first experiment, we have calculated the logarithm of the following non-singular triangular matrices, with no real negative eigenvalues (MATLAB style is used) by the MATLAB function \texttt{logm}, without and with the scaling technique:

\begin{itemize}
	\item[$\bullet$] $T_1$ is the matrix $\mathtt{exp(a)*[1\ b;0\ 1]}$, with $a=0.1$ and $b=10^6$; its exact logarithm is $\mathtt{[a\ b;0\ a]}$;
	\item[$\bullet$] $T_2$ is the matrix in (\ref{T1});
	\item[$\bullet$] $T_3$ has been obtained by $\mathtt{[\sim,T3]=schur(gallery('frank',8),'complex')}$;
	\item[$\bullet$] $T_4$ has been obtained by $\mathtt{[\sim,T4]=schur(gallery('dramadah',11),'complex')}$;
	\item[$\bullet$] $T_5$ has come from $\mathtt{[\sim,T5]=schur(gallery('frank',13),'complex')}$;
	\item[$\bullet$] $T_6$ to $T_{10}$ are randomized matrices with orders ranging from $9$ to $20$ with smaller entries in the diagonal than in the super-diagonals.
\end{itemize}

\begin{figure}
	\centering
	\includegraphics[width=15cm]{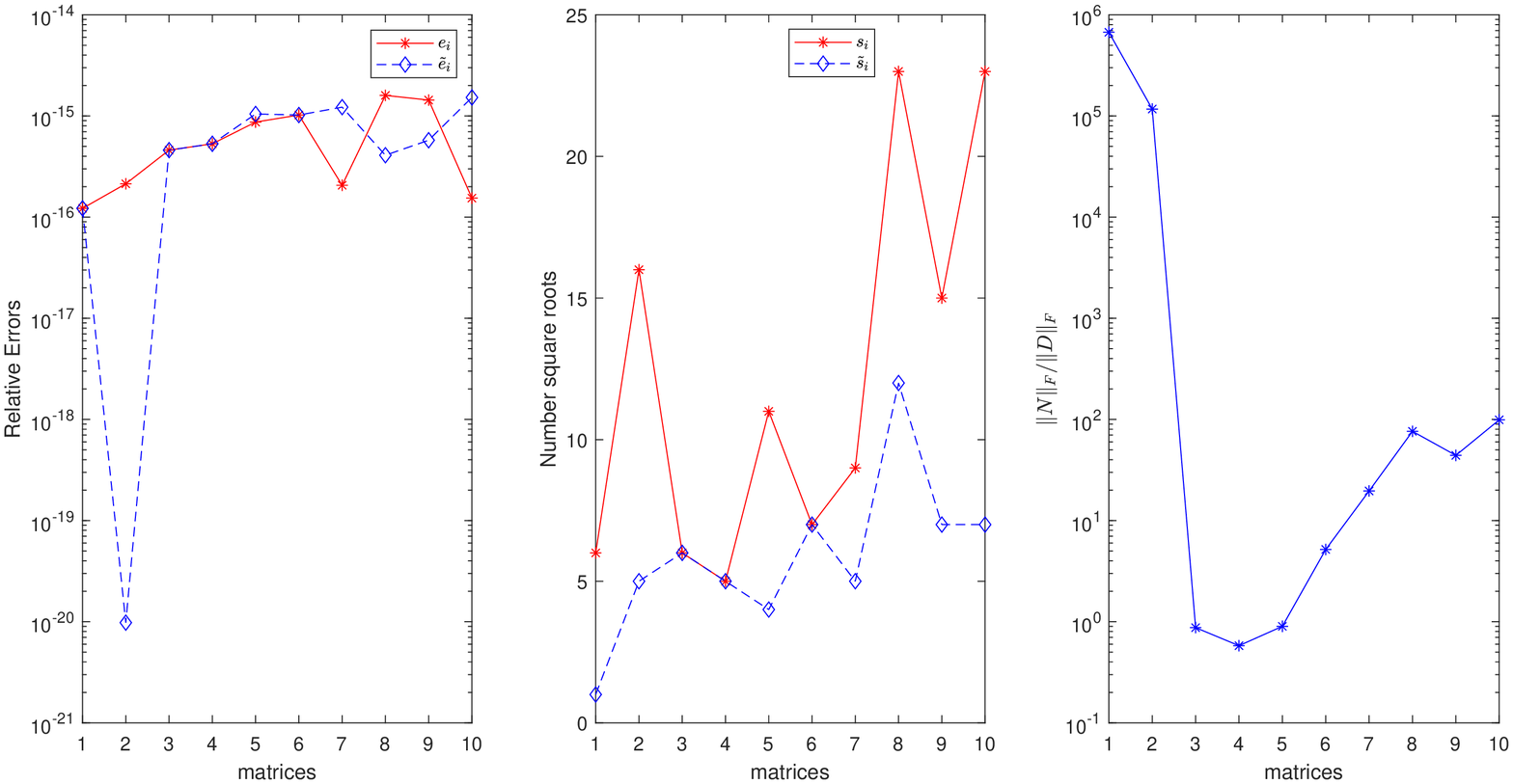}
	\caption{\small Results of Experiment 1. Left: Relative errors of the approximations obtained by \texttt{logm} with the scaling technique ($\tilde{e}$) and without  ($e$). Middle: Number of square roots required by \texttt{logm} with the scaling technique ($\tilde{s}$) and without ($s$). Right: Values of the ratio $\|N\|_F/\|D\|_F$.}
	\label{fig-log}
\end{figure}

The results are displayed in Figure \ref{fig-log}. Excepting the matrix $T_1$, whose exact logarithm is known, we have considered as the exact $\log(T_i)$ ($i=2,\ldots,10$) the result of evaluating the logarithm at 200 decimal digit precision using the Symbolic Math Toolbox and rounding the result to double precision.

\begin{figure}
	\centering
	\includegraphics[width=15cm]{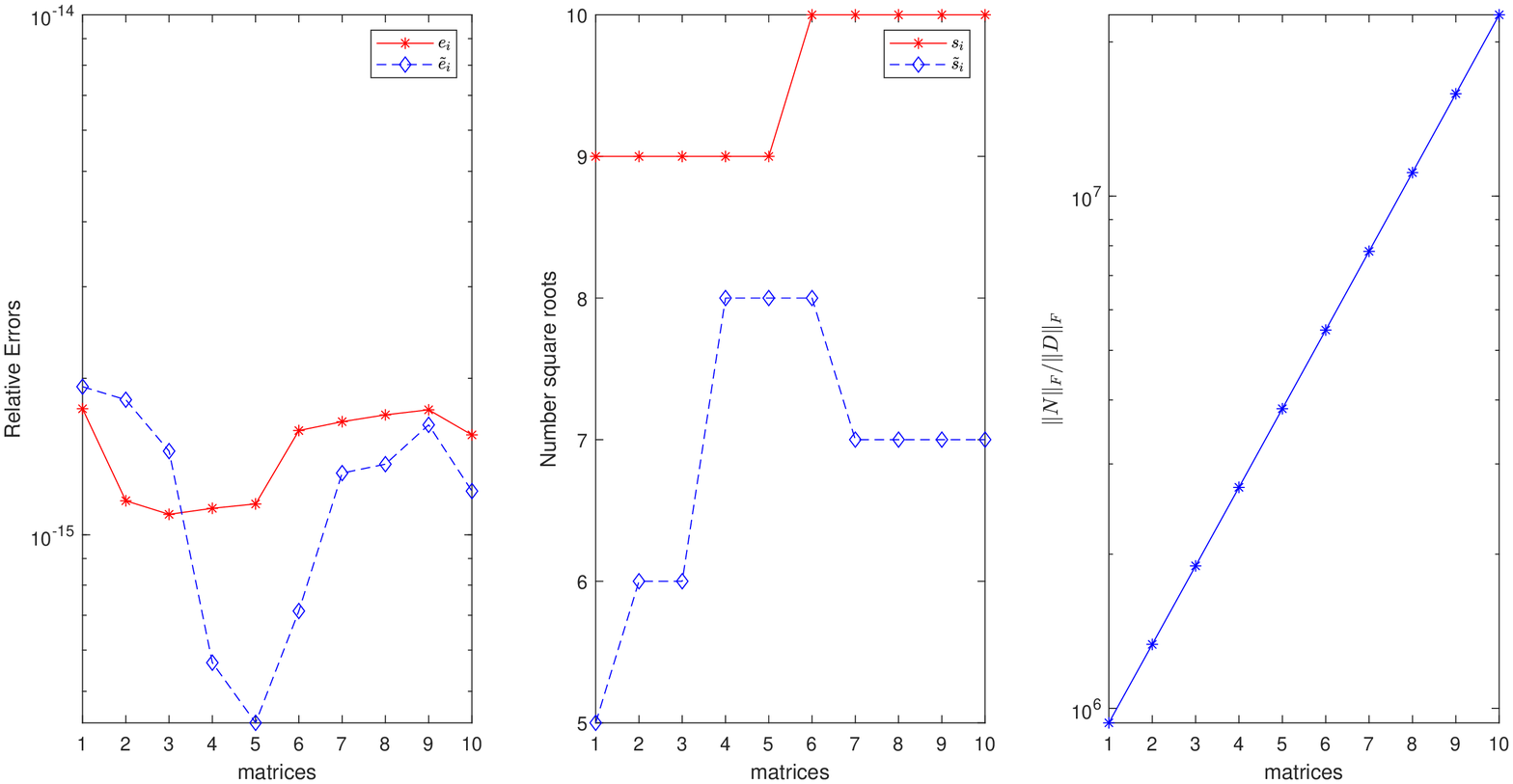}
	\caption{\small Results of Experiment 2. Left: Relative errors of the approximations obtained by \texttt{logm} with the scaling technique ($\tilde{e}$) and without  ($e$). Middle: Number of square roots required by \texttt{logm} with the scaling technique ($\tilde{s}$) and without ($s$). Right: Values of the ratio $\|N\|_F/\|D\|_F$.}
	\label{fig-log-toep}
\end{figure}

\medskip\noindent {\it Experiment 2.} This experiment illustrates the performance of Algorithm \ref{alg1} for the logarithm of matrices having larger size. We consider ten Toeplitz matrices $T=[t_{ij}]$ ($i,j=1,\ldots,n$)  with sizes $n=82:2:100$. Its entries are defined by
$$t_{ij}=\left\{\begin{array}{lll}
	1.2^{j-i+1} &\mbox{if} & i\leq j\\
	0 & \mbox{if} & i> j
	\end{array} \right.
$$
and the results are displayed in Figure \ref{fig-log-toep}.	

\begin{figure}
 	\centering
 	\includegraphics[width=15cm]{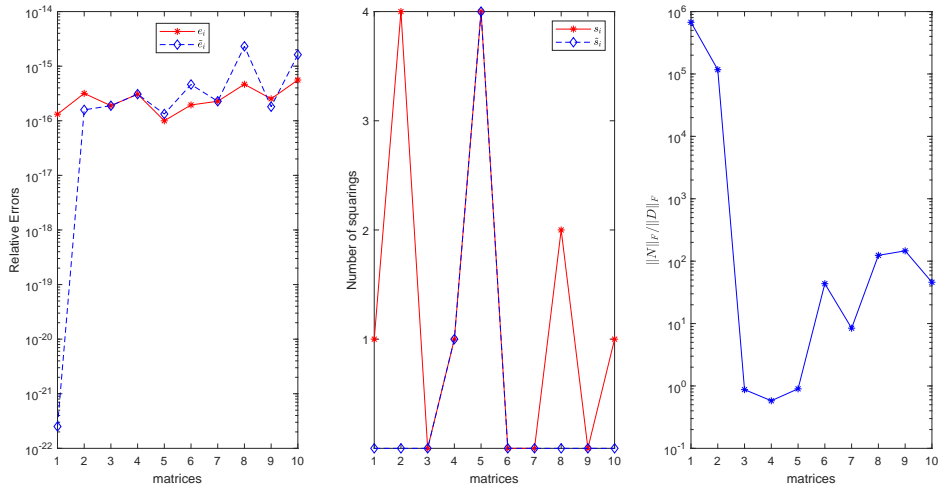}
 	\caption{\small Results of Experiment 3. Left: Relative errors of the approximations obtained by \texttt{expm} with the scaling technique ($\tilde{e}$) and without  ($e$). Middle: Number of squarings required by \texttt{expm} with the scaling technique ($\tilde{s}$) and without ($s$). Right: Values of the ratio $\|N\|_F/\|D\|_F$.}
 	\label{fig-exp}
 \end{figure}

\medskip\noindent {\it Experiment 3.} The results of this experiment are displayed in Figure \ref{fig-exp} and concern to the computation of the exponential of the same matrices of Experiment 1, with and without the scaling technique.

\begin{figure}
	\centering
	\includegraphics[width=15cm]{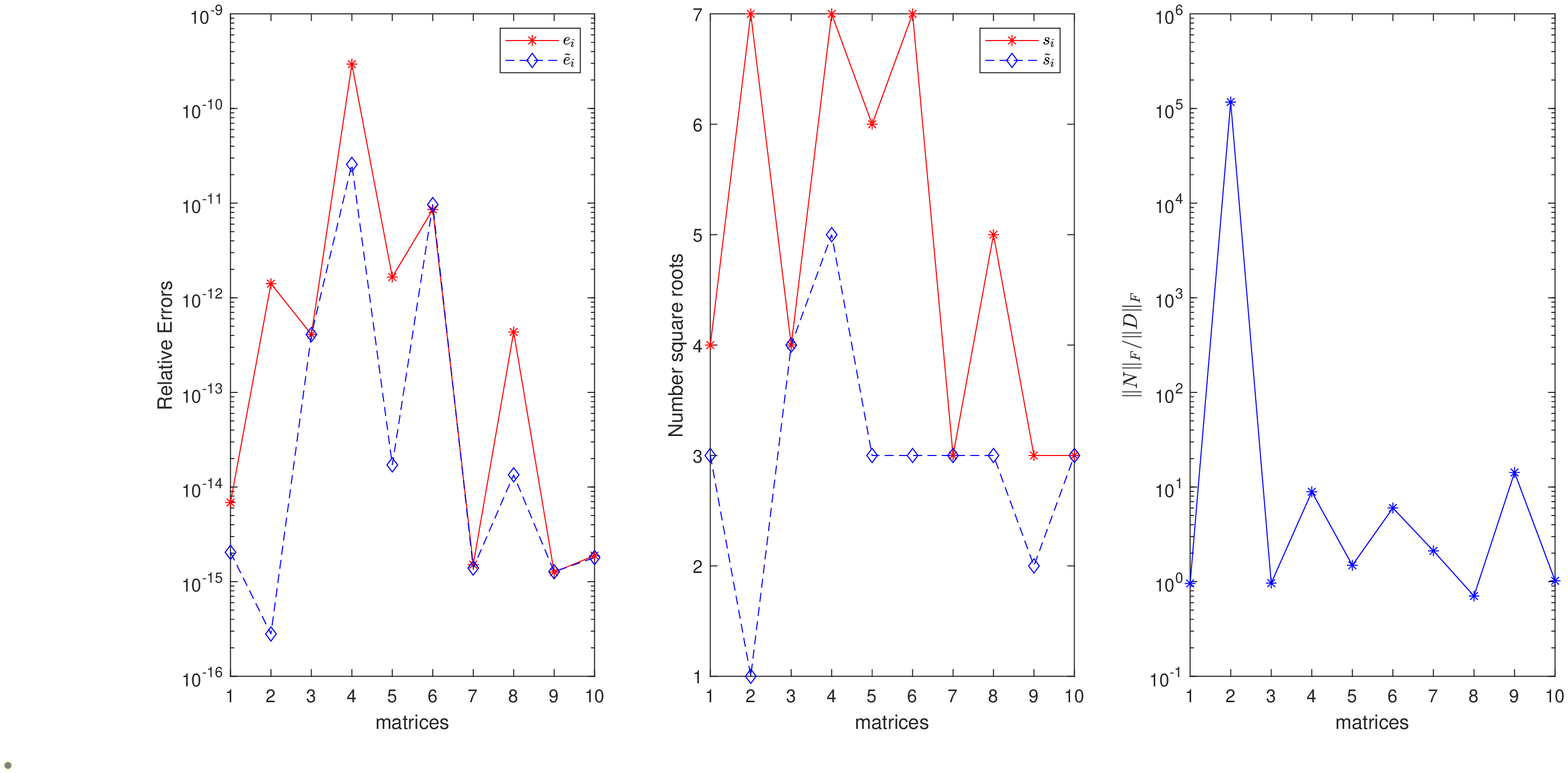}
	\caption{\small Results of Experiment 4. Left: Relative errors of the approximations obtained by \texttt{acosm} with scaling ($\tilde{e}$) and without ($e$). Middle: Number of square roots required by \texttt{acosm} with the scaling technique ($\tilde{s}$) and without ($s$). Right: Values of the ratio $\|N\|_F/\|D\|_F$.}
	\label{fig-acos}
\end{figure}

\medskip\noindent {\it Experiment 4.} In this experiment, we have implemented Algorithm 5.2 of \cite{Aprahamian16} for computing the inverse cosine of a matrix, $\acos(T)$, with and without scaling. MATLAB codes are available at https://github.com/higham/matrix-inv-trig-hyp and the algorithm is called by \texttt{acosm}. It involves a prior Schur decomposition, thus being well suited to be combined with the scaling technique. We have considered ten matrices from MATLAB's gallery with no eigenvalues in the set $\{-1,\,1\}$ and sizes ranging from $30$ to $50$. The results can be observed in Figure \ref{fig-acos}.

\medskip The results of Experiment 1 and 2, which are displayed in Figures \ref{fig-log} and \ref{fig-log-toep}, for the MATLAB's function \texttt{logm}, show that the scaling technique may bring a reduction on the number of square roots, especially in upper triangular matrices with small entries in diagonal if compared with the magnitude of the entries in super-diagonals. Regarding the relative errors, we may observe slight fluctuations, but in general the relative errors are not affected significantly by the scaling technique. However, if instead of \texttt{logm}, which involves sophisticated techniques \cite{Mohy09}, we combine the proposed scaling technique with the standard inverse scaling and squaring of \cite{Kenney}, then we could see a reduction not only on the number of square roots but also on relative errors.

In Experiment 3 for \texttt{expm} (Figure \ref{fig-exp}), the effects of the scaling technique are not so observable as with the \texttt{logm}. We recall that the matrices are the same as in Experiment 1.

In Experiment 4 for \texttt{acosm} (see Figure \ref{fig-acos}), we see a slight improvement on relative errors and less square roots in many cases. Almost the experiments suggest the existence of a connection between the savings in the number of square roots or squarings and the value of the ratio $\|N\|_F/\|D\|_F$.

\section{Conclusions}\label{conclusions}
We have proposed an simple and inexpensive scaling technique aimed at improving algorithms for evaluating functions of triangular matrices, in terms of computational cost and accuracy. It is particularly well suited to be combined with algorithms involving a prior Schur decomposition. Such a technique involves a scalar $\alpha$ that needs to be carefully chosen. We have presented a practical strategy for finding such a scalar, that has given promising results for experiments involving the matrix exponential, the logarithm and the inverse cosine. However, further research is needed for tuning appropriate values for $\alpha$ and $m$ and for understanding the effects of the proposed technique on algorithms for computing other matrix functions (check the list of algorithms and functions mentioned at the end of Section \ref{pre}).

\end{document}